\documentclass[11pt,reqno]{amsart}

\usepackage{geometry}
\usepackage{amsmath,amsthm,amssymb,amsfonts,dsfont}
\usepackage{mathtools}
\usepackage{float}
\usepackage[utf8]{inputenc} 
\usepackage[T1]{fontenc}    
\usepackage[hidelinks]{hyperref}
\usepackage[capitalise]{cleveref}

\usepackage{comment}
\usepackage{subcaption}
\usepackage{wrapfig}
\usepackage{bm}
\usepackage{url}           
\usepackage{booktabs}    
\usepackage{nicefrac}    
\usepackage{microtype}    
\usepackage{bbm}
\usepackage{cite}
\usepackage{psfrag}
\usepackage{cite}
\usepackage{color,soul}
\usepackage{breakcites}   
\usepackage{bigints}
\usepackage{enumitem}
\setlist{leftmargin=1.6em}

\newtheorem{remark}{Remark}
\newtheorem{theorem}{Theorem}
\newtheorem{proposition}{Proposition}
\newtheorem{lemma}{Lemma} 

\usepackage{graphicx}
\usepackage{xcolor}
\usepackage{svg}
\usepackage{algorithm}
\usepackage[noend]{algpseudocode}

\makeatletter
\def\subsubsection{\@startsection{subsubsection}{3}%
  \z@{.5\linespacing\@plus.7\linespacing}{-.5em}%
  {\normalfont\bfseries}}
\makeatother

\newcommand{\vertiii}[1]{{\left\vert\kern-0.25ex\left\vert\kern-0.25ex\left\vert #1 
    \right\vert\kern-0.25ex\right\vert\kern-0.25ex\right\vert}}

\newcommand{\sH}{\mathsf{H}}

\newcommand{\cC}{\mathcal{C}}

\newcommand{\cF}{\mathcal{F}}

\newcommand{\cP}{\mathcal{P}}

\newcommand{\cU}{\mathcal{U}}

\newcommand{\cX}{\mathcal{X}}
\newcommand{\cY}{\mathcal{Y}}

\newcommand{\RR}{\mathbb{R}}

\newcommand{\vasti}{\bBigg@{3.5 }}
\newcommand{\vast}{\bBigg@{4}}
\newcommand{\Vast}{\bBigg@{5}}
\newcommand{\Vastt}{\bBigg@{7}}

\makeatother

\newcommand{\be}{\begin{equation}}
\newcommand{\ee}{\end{equation}}
\newcommand{\ba}{\begin{align}}
\newcommand{\ea}{\end{align}}
\newcommand{\baa}{\begin{align*}}
\newcommand{\eaa}{\end{align*}}

\newcommand{\KL}{\mathsf{D}_{\mathsf{KL}}}
\newcommand{\OT}{\mathsf{OT}}

\usepackage{geometry}
\usepackage{mathtools}

\usepackage{amsmath,amsthm,amssymb,amsfonts,bbm,dsfont}
\usepackage{float}

\usepackage{comment}
\usepackage{subcaption}
\usepackage{bm}
\usepackage{color,soul}
\usepackage{breakcites}     
\usepackage{bigints}

\usepackage{enumitem}

\usepackage{graphicx}
\usepackage{xcolor}
\usepackage{svg}
\usepackage{algorithm}
\usepackage[noend]{algpseudocode}
\usepackage{caption}
\usepackage{subcaption}

\makeatletter

\newcommand{\eps}{\varepsilon}

\makeatother

\begin{document}

\title[Neural Estimation of Entropic Optimal Transport]{Neural Estimation of Entropic Optimal Transport}

\date{First version: \today}

\author[T. Wang]{Tao Wang}

\address[T. Wang]{Applied Mathematics and Computational Science, University of Pennsylvania.}
\email{tawan@sas.upenn.edu}

\author[Z. Goldfeld]{Ziv Goldfeld}
\address[Z. Goldfeld]{School of Electrical and Computer Engineering, Cornell University.}
\email{goldfeld@cornell.edu}

\begin{abstract}
Optimal transport (OT) serves as a natural framework for comparing probability measures, with applications in statistics, machine learning, and applied mathematics. Alas, statistical estimation and exact computation of the OT distances suffer from the curse of dimensionality. To circumvent these issues, entropic regularization has emerged as a remedy that enables parametric estimation rates via plug-in and efficient computation using Sinkhorn iterations. Motivated by further scaling up entropic OT (EOT) to data dimensions and sample sizes that appear in modern machine learning applications, we propose a novel neural estimation approach. Our estimator parametrizes a semi-dual representation of the EOT distance by a neural network, approximates expectations by sample means, and optimizes the resulting empirical objective over parameter space. We establish non-asymptotic error bounds on the EOT neural estimator of the cost and optimal plan. Our bounds characterize the effective error in terms of neural network size and the number of samples, revealing optimal scaling laws that guarantee parametric convergence. The bounds hold for compactly supported distributions, and imply that the proposed estimator is minimax-rate optimal over that class. Numerical experiments validating our theory are also provided.
\end{abstract}

\thanks{Z. Goldfeld is partially supported by NSF grants CCF-2046018, DMS-2210368, and CCF-2308446, and the IBM Academic Award.
}

\maketitle
\section{Introduction}
Optimal transport (OT) theory \cite{villani2009optimal} provides a natural framework for comparing probability distributions. Specifically, given two Borel probability measures $\mu, \nu$ on $\mathbb{R}^d$, the OT problem between them with cost function $c$ is defined as
\begin{equation}
    \label{eq:OT}
\OT_c(\mu, \nu):=\inf _{\pi \in \Pi(\mu, \nu)} \int_{\mathbb{R}^d \times \mathbb{R}^d} c(x, y) d \pi(x, y)
\end{equation}
where $\Pi(\mu, \nu)$ is the set of couplings between $\mu$ and $\nu$. The special case is the $p$-Wasserstein distance for $p \in[1, \infty)$, is given by $\mathsf{W}_p(\mu, \nu):=\left(\OT_{\|\cdot\|^p}(\mu, \nu)\right)^{1 / p}$. 
The Wasserstein distance has found applications in various fields, encompassing machine learning \cite{arjovsky2017wasserstein,tolstikhin2017wasserstein,courty2017joint}, statistics \cite{carlier2016vector,chernozhukov2017monge,ghosal2022multivariate}, and applied
mathematics \cite{jordan1998variational,santambrogio2017euclidean}. 
This widespread applicability is driven by an array of desirable properties that the Wasserstein distance possesses, including its metric structure ($\mathsf{W}_p$ metrizes weak convergence plus convergence of $p$-th moments), a convenient dual form, robustness to support mismatch, and a rich geometry it induces on a space of probability measures.

Despite the aforementioned empirical progress, the OT problem suffers from the statistical and computational hardness issues. The estimation rate of the OT cost between distributions on $\RR^d$ is generally $n^{-1/d}$ (without further assumptions) \cite{fournier2015rate}, which deteriorates exponentially with dimensions---a phenomenon known as the curse of dimensionality. Computationally, OT is a linear program (LP), solvable in $O(n^3 \log (n))$ time for distribution on $n$ points using interior point methods or min cost flow algorithms \cite{peyre2017computational}. However, as statistical considerations mandate $n$ to scale exponentially with $d$ to get accurate estimates, the LP computational paradigm becomes infeasible when dimension is large. To circumvent these issues, entropic regularization has emerged as a popular alternative \cite{cuturi2013sinkhorn}
\begin{equation}
    \label{eq:EOT}
    \OT_{c}^ \varepsilon(\mu, \nu):=\inf _{\pi \in \Pi(\mu, \nu)} \int c\, d \pi+\varepsilon\KL (\pi \| \mu \otimes \nu),
\end{equation}
where $\KL$ is the Kullback-Leibler divergence and $\varepsilon>0$ is a regularization parameter. Empirical estimation of EOT enjoys the parametric $n^{-1/2}$ convergence rate in arbitrary dimension, under several settings \cite{genevay2019sample,mena2019statistical}. Computationally, EOT between discrete distributions can be efficiently solved via the Sinkhorn algorithm \cite{cuturi2013sinkhorn} in $O(n^2)$ time. However, even this quadratic time complexity is prohibitive when dealing with large and high-dimensional datasets that appear in modern machine learning tasks. Motivated to scale up EOT to such regimes, this work develops a novel neural estimation approach that is end-to-end trainable via backpropagation, compatible with minibatch-based optimization, and adheres to strong performance guarantees. 

\subsection{Contributions}

We focus on the nominal case of the quadratic EOT distance, i.e., $c(x,y)=\frac{1}{2}\|x-y\|^2$. Thanks to the EOT semi-dual form, we have
\begin{equation}
    \label{eq:semi_dual}
    \OT_{c}^\varepsilon(\mu,\nu)=\sup_{\varphi\in L^1(\mu)} \int_{\RR^d}\varphi\, d\mu+\int_{\RR^d}\varphi^{c, \varepsilon}d\,\nu,
\end{equation}
where $\varphi^{c,\varepsilon}$ is $(c,\varepsilon)$-transform of $\varphi$  with respect to (w.r.t.) the cost function. We study regularity of optimal dual potentials $\varphi$ and show that they belong to a H\"older class of arbitrary smoothness. Leveraging this, we define our neural estimator (NE) by parametrizing the dual potential using a neural network (NN), approximating expectations by sample means, and optimizing the resulting empirical objective over the NN parameters. Our approach yields not only an estimate of the EOT distance, but also a neural EOT plan that is induced by the learned NN. As the estimator is trainable via gradient methods using backpropagation and minibatches, it can seamlessly integrated into downstream tasks as a loss, a regularizer, or a discrepancy quantification module.

We provide formal guarantees on the quality of the NE of the EOT cost and the corresponding transportation plan. Our analysis relies on non-asymptotic function approximation theorems and tools from empirical process theory to bound the two sources of error involved: function approximation and empirical estimation. Given $n$ samples from the population distributions, we show that the effective error of a NE realized by a shallow NN of $k$ neurons scales as
\begin{equation}
 \label{eq:stat_rate}   
O\left(\operatorname{poly}(1 / \varepsilon)\left(k^{-1 / 2}+n^{-1 / 2}\right)\right)
\end{equation}
with the polynomial dependence on $1 / \varepsilon$ explicitly characterized. This bound on the EOT cost estimation error holds for arbitrary, compactly supported distributions. This stands in struck contrast to existing neural estimation error bounds for other divergences \cite{nguyen2010estimating,sreekumar2021non,sreekumar2022neural,tsur2023max}, which typically require strong regularity assumptions on the population distributions (e.g., H\"older smoothness of densities). This is unnecessary in our setting thanks to the inherit regularity of dual EOT potentials for smooth cost functions, such as our quadratic cost.

The above bound reveals the optimal scaling of the NN and dataset sizes, namely $k \asymp$ $n$, which achieves the parametric convergence rate of $n^{-1 / 2}$ and guarantees minimax-rate optimality of our NE. The explicitly characterized polynomial dependence on $\varepsilon$ in our bound matches the bounds for EOT estimation via empirical plug-in \cite{mena2019statistical,groppe2023lower}. We also note that our neural estimation results readily extend to the EOT problem with general smooth cost functions. The developed NE is empirically tested on synthetic data, demonstrating its scalability to high dimensions and validating our theory.

\subsection{Related Literature}\label{subsec:lietrature_review}

Neural estimation is a popular approach for enhancing scalability. Prior research explored the tradeoffs between approximation and estimation errors in non-parametric regression \cite{barron1994approximation,bach2017breaking,suzuki2018adaptivity} and density estimation \cite{yang1999information,uppal2019nonparametric} tasks. More recently, neural estimation of statistical divergences and information measures has been gaining attention. The mutual information NE (MINE) was proposed in \cite{belghazi2018mine}, and has seen various improvements since \cite{poole2018variational,song2019understanding,chan2019neural,mroueh2021improved}. Extensions of the neural estimation approach to directed information were studied in \cite{molavipour2021neural,tsur2023neural,tsur2023data}. Theoretical guarantees for $f$-divergence NEs, accounting for approximation and estimation errors, as we do here, were developed in \cite{sreekumar2021non,sreekumar2022neural} (see also \cite{nguyen2010estimating} for a related approach based on reproducing kernel Hilbert space parameterization). Neural estimation of the Stein discrepancy and the minimum mean squared error were considered in \cite{repasky2023neural} and \cite{diaz2021lower}, respectively. Neural methods for approximate computation  of the Wasserstein distances have been considered under the Wasserstein generative adversarial network (GAN) framework \cite{arjovsky2017wasserstein,gulrajani2017improved}, although these approaches are heuristic and lack formal guarantees. Utilizing entropic regularization, \cite{daniels2021score} studied a score-based generative neural EOT model, while an energy-based model was considered in \cite{mokrov2023energy}.

\section{Background and Preliminaries}

\subsection{Notation}

Let $\|\cdot\|$ and $\langle \cdot, \cdot \rangle$ designate Euclidean norm and the inner product in $\RR^d$, respectively. For $1 \leq p < \infty$, the $L^p$ space over $\mathcal{X} \subseteq \mathbb{R}^d$ with respect to (w.r.t.) the measure $\mu$ is denoted by $L^p(\mu)$, with $\|f\|_{p, \mu}\coloneqq\big( \int_{\cX}|f|^pd\mu\big)^{1/p}$ representing the norm. We use $\|\cdot\|_{\infty, \mathcal{X}}$ for standard sup-norm on $\mathcal{X} \subseteq \mathbb{R}^d$ (i.e., when  $p=\infty$). Slightly abusing notation, we also set $\|\mathcal{X}\|\coloneqq\sup _{x \in \mathcal{X}}\|x\|_{\infty}$. The class of Borel probability measures on $\mathcal{X} \subseteq \mathbb{R}^d$ is denoted by $\mathcal{P}(\mathcal{X})$. For $\mu, \nu \in \mathcal{P}(\mathcal{X})$ with $\mu \ll \nu$, i.e., $\mu$ is absolutely continuous w.r.t. $\nu$, we use $\frac{d\mu}{d\nu}$ for the Radon-Nikodym derivative of $\mu$ w.r.t. $\nu$. The subset of probability measures that are absolutely continuous w.r.t. the Lebesgue measure is denoted by $\mathcal{P}_{\mathsf{ac}}(\cX)$.
We use $\lesssim_x$ to denote inequalities up to constants that only depend on $x$; the subscript is dropped when the constant is universal. For $a,b \in \mathbb{R}$, we write $a \vee b = \max \{ a,b \}$ and $a \land b = \min \{ a,b \}$.

Some additional notation used for our derivations are as follows. For any multi-index $\alpha=(\alpha_1,\dots,\alpha_d) \in \mathbb{N}_0^d$ with $|\alpha|=\sum_{j=1}^d \alpha_j\left(\mathbb{N}_0=\mathbb{N} \cup\{0\}\right)$, define the differential operator $D^\alpha = \frac{\partial^{|\alpha|}}{\partial x_1^{\alpha_1} \cdots \partial x_{d}^{\alpha{}_d}}$
with $D^0 f = f$. We write $N(\delta, \cF, \mathsf{d})$ for the $\delta$-covering number of a function class $\cF$ w.r.t. a metric $\mathsf{d}$, and $N_{[\,]}(\delta, \cF, \mathsf{d})$ for the bracketing number. For an open set $\mathcal{U} \subseteq \mathbb{R}^d$, $b \geq 0$, and an integer $m \geq 0$, let $\cC_b^m(\mathcal{U})\coloneqq \big\{f \in \cC^m(\mathcal{U}):\,\max _{\alpha:|\alpha| \leq m}\left\|D^\alpha f\right\|_{\infty, \mathcal{U}} \leq b\big\}$ denote the H\"older space of smoothness index $m$ and radius $b$. The restriction of $f: \mathbb{R}^d \rightarrow \mathbb{R}$ to a subset $\mathcal{X} \subseteq \mathbb{R}^d$ is denoted by $f\big|_{\mathcal{X}}$.

\subsection{Entropic Optimal Transport}
We briefly review basic definitions and results concerning EOT problems. Let $\cX \subseteq \RR^d$, given distributions $(\mu,\nu)\in\cP(\cX)\times \cP(\cY)$ and a cost function $c:\cX\times \cY\to\RR$, the primal EOT formulation is obtained by regularizing the OT cost by the KL divergence, 
\begin{equation}
\label{EOT}
      \OT_c^{\eps}(\mu,\nu)\coloneqq \inf_{\pi\in\Pi(\mu,\nu)} \int_{\cX\times\cY} c\,d\pi+\eps \KL(\pi\|\mu\otimes \nu),
\end{equation}
where $\eps>0$ is a regularization parameter and $\KL(\mu\|\nu)\coloneqq\int \log \left(\frac{d \mu}{d \nu}\right) d \mu$ if $\mu \ll \nu$ and $+\infty$ otherwise. Classical OT \cite{villani2009optimal} is obtained from \eqref{EOT} by setting $\eps=0$. When~$c\in L^1(\mu\otimes \nu)$, EOT admits the dual and semi-dual formulations, which are, respectively, given~by
\begin{align}
\OT_c^{\eps}(\mu,\nu)&=\sup_{(\varphi,\psi)\in L^1(\mu)\times L^1(\nu)}\int \varphi d\mu+\int \psi d\nu-\eps \int e^{\frac{\varphi\oplus \psi-c}{\eps}}d\mu\otimes \nu+\eps, \label{eq:EOT_dual}\\
&=\sup_{\varphi\in L^1(\mu)}\int \varphi d\mu+\int \varphi^{c,\eps}d\nu,\label{eq:EOT_semidual}
\end{align}
where we have defined $(\varphi\oplus \psi)(x,y)=\varphi(x)+\psi(y)$ and the $(c,\eps)$-transform of $\varphi$ is given by $\varphi^{c,\eps}=-\eps \log \left(\int_{\cX} \exp \left(\frac{\varphi(x)-c(x,\cdot)}{\eps}\right) d \mu(x)\right)$. There exist functions $(\varphi,\psi)$ that achieve the supremum in $\eqref{eq:EOT_dual}$, which we call \emph{EOT potentials}. These potentials are almost surely (a.s.) unique up to additive constants, i.e., if $(\tilde{\varphi}, \tilde{\psi})$ is another pair of EOT potentials, then there exists a constant $a \in \mathbb{R}$ such that $\tilde{\varphi}=\varphi+a$ $\mu$-a.s. and $\tilde{\psi}=\psi-a$ $\nu$-a.s. 

A pair $(\varphi, \psi) \in L^1(\mu) \times L^1(\nu)$ are EOT potentials if and only if they satisfy the Schr\"odinger system
\begin{equation}
    \label{eq:Schrodinger system}
\int e^{\frac{\varphi(x)+\psi(\cdot)-c(x, \cdot)}{\eps}} d \mu(x)=1 \quad \nu \text {-a.s. } \quad \text { and } \quad \int e^{\frac{\varphi(\cdot)+\psi(y)-c(\cdot, y)}{\eps}} d \nu(y)=1 \quad \mu \text {-a.s. }
\end{equation}
Furthermore, $\varphi$ solves the semi-dual from \eqref{eq:EOT_semidual} if an only if $(\varphi,\varphi^{c,\eps})$ is a solution to the full dual in \eqref{eq:EOT_dual}. Given EOT potentials $(\varphi, \psi)$, the unique EOT plan can be expressed in their terms as $d \pi^\eps_{\star}=e^{\frac{\varphi\oplus \psi-c}{\eps}} d \mu \otimes \nu$. Subject to smoothness assumptions on the cost function and the population distributions, various regularity properties of EOT potentials can be derived; cf., e.g., \cite[Lemma 1]{goldfeld2022limit}.    

\section{Neural Estimation of EOT Cost and Plan}
We consider compactly supported distributions $(\mu, \nu) \in \mathcal{P}(\mathcal{X}) \times \mathcal{P}(\mathcal{Y})$ and the quadratic cost function $c(x,y)=\frac{1}{2}\|x-y\|^2$ (henceforth dropping the subscript $c$). For simplicity, further assume that $\cX,\cY\subseteq[-1,1]^{d}$, although our results readily extend to arbitrary compact supports. We next describe the NE for the EOT distance and plan, followed by non-asymptotic performance guarantees for both objects. All proofs are deferred to the supplement.

\subsection{EOT Neural Estimator}

For $(\mu, \nu) \in \mathcal{P}(\cX) \times \mathcal{P}(\cY)$, let $X^n\coloneqq\left(X_1, \cdots, X_n\right)$ and $Y^n\coloneqq\left(Y_1, \cdots, Y_n\right)$ be $n$ independently and identically distributed (i.i.d.) samples from $\mu$ and $\nu$, respectively. Further suppose that the sample sets are independent of each other. Denote the empirical measures induced by these samples as $\hat\mu_n=n^{-1}\sum_{i=1}^n\delta_{X_i}$ and $\hat\nu_n=n^{-1}\sum_{i=1}^n\delta_{Y_i}$.

Our NE is realized by a shallow ReLU NN (i.e., a single hidden layer) with $k$ neurons, which defines the function class
\begin{equation}
\label{NN_class}
    \cF_{k,d}(a)\coloneqq \left\{f: \mathbb{R}^{d} \rightarrow \mathbb{R}:\begin{aligned}
&f(x)=\sum_{i=1}^k \beta_i \phi\left(w_i \cdot x+b_i\right)+w_0 \cdot x+b_0, \\
& \max _{1 \leq i \leq k}\left\|w_i\right\|_1 \vee\left|b_i\right| \leq 1,\ \max _{1 \leq i \leq k}\left|\beta_i\right| \leq 2ak^{-1},\ \left|b_0\right| \leq a,\ \left\|w_0\right\|_1 \leq a
\end{aligned}\right\},
\end{equation}
where $a\in\RR_{\geq 0}$ specifies the parameter bounds and $\phi:\RR\to\RR_{\geq 0}: z\mapsto z\vee 0$ is the ReLU activation function, which acts on vectors component-wise.

We parametrize the semi-dual form of $\OT^{\mspace{1mu}\eps}(\mu,\nu)$ (see \eqref{eq:EOT_semidual}) using a NN from the class $\cF_{k,d}(a)$ and replace expectations with sample means. Specifically, the EOT distance NE is
\begin{equation}\label{eq:NE_EOT}
\begin{aligned}
    &\widehat\OT_{ k, a}^{\mspace{1mu}\eps}(X^n, Y^n)\coloneqq 
   \max_{f \in \cF_{k, d}(a)} \frac{1}{n}\sum_{i=1}^nf(X_i)-\frac{\eps}{n}\sum_{j=1}^n\log\left(\frac{1}{n}\sum_{i=1}^n\exp\left(\frac{f(X_i)-\frac{1}{2}\|X_i-Y_j\|^2}{\eps}\right)\right). 
\end{aligned}
\end{equation}
For any NN $f\in \cF_{k,d}(a)$, we define the induced neural plan
\begin{equation}
    d \pi^\eps_{f}(x, y) \coloneqq \frac{\exp \left(\frac{f(x)-\frac{1}{2}\|x-y\|^2}{\eps}\right)}{\int_{\cX} \exp \left(\frac{f(x)-\frac{1}{2}\|x-y\|^2}{\eps}\right) d\mu(x)} d\mu\otimes\nu(x,y).\label{eq:NE_plan}
\end{equation}
Upon computing the NE in \eqref{eq:NE_EOT}, the neural plan $d \pi^\eps_{f_\star}$ induced by an optimal NN $f_\star\in\cF_{k,d}(a)$ serves as an estimate of the true optimal plan $\pi^\eps_\star$ that achieves the infimum in \eqref{EOT}.

\subsection{Performance Guarantees}

We provide formal guarantees for the neural estimator of the EGW cost and the neural transportation plan defined above. Starting from the cost estimation setting, we establish two separate bounds on the effective (approximation plus estimation) error. The first is non-asymptotic and presents optimal convergence rates, but calibrates the NN parameters to a cumbersome dimension-dependent constant. Following that, we present an alternative bound that avoids the dependence on the implicit constant, at the expense of a polylogarithmic slow-down in the rate and a requirement that the NN size $k$ is large enough. 

\begin{theorem}[EOT cost neural estimation; bound 1] \label{thm:Bound_EGW1}
     There exists a constant $C>0$ depending only on $d$, such that setting $a=C(1+\eps^{1-s})$ with $s=\left\lfloor d  / 2\right\rfloor+3$, we have
\begin{equation}
\label{Bound_EGW1}
\begin{aligned}
&\sup _{(\mu, \nu) \in \mathcal{P}(\cX) \times \mathcal{P}(\cY)} \mathbb{E}\left[\left|\widehat\OT_{k,a}^{\mspace{1mu}\eps}(X^n, Y^n)-\OT^{\eps}(\mu, \nu)\right|\right]\\
&\quad\quad\quad\quad\quad\quad\lesssim_{d} \left(1+\frac{1}{\eps^{\left\lfloor \frac{d}{2}\right\rfloor+2}}\right) k^{-\frac{1}{2}}+\min\left\{1+\frac{1}{\eps^{\left\lceil \frac{3d}{2}\right\rceil+4}}\,,\left(1+\frac{1}{\eps^{\left\lfloor \frac{d}{2}\right\rfloor+2}}\right)\sqrt{k}\right\} n^{-\frac{1}{2}}.
\end{aligned}
\end{equation}
\end{theorem}

The proof of Theorem \ref{thm:Bound_EGW1} is given in Appendix \ref{appen:proof_Bound_EGW1}. We establish  regularity of semi-dual EOT potentials (namely, $(\varphi,\varphi^{c,\eps})$ in \eqref{eq:EOT_semidual}), showing that they belong to a H\"older class of arbitrary smoothness. This, in turn, allows accurately approximating these dual potentials by NNs from the class $\cF_{k,d}(a)$ with error $O(k^{-1/2})$, yielding the first term in the bound. To control the estimation error, we employ standard maximal inequalities from empirical process theory along with a bound on the covering or bracketing number of $(c,\eps)$-transform of the NN class. The resulting empirical estimation error bound comprises the second term on the RHS above.

\begin{remark}[Almost explicit expression for $C$]
The expression of the constant $C$ in Theorem \ref{thm:Bound_EGW1} is cumbersome, but can nonetheless be evaluated. Indeed, one may express $C=C_{s} C_{d} \bar{c}_{d}$, with explicit expressions for $C_{d}$ and $\bar{c}_{d}$ given in \eqref{eq:constant_C_dx_dy} and \eqref{eq:constant_c_bar}, respectively, while $C_s$ is a combinatorial constant that arises from the multivariate Faa di Bruno formula (cf. \eqref{eq:constant_C_s_1}-\eqref{eq:constant_C_s_2}). The latter constant is quite convoluted and is the main reason we view $C$ as implicit. 
\end{remark}

Our next bound circumvents the dependence on $C$ by letting the NN parameters grow with its size $k$. This bound, however, requires $k$ to be large enough and entails additional polylog factors in the rate. The proof is similar to that of Theorem \ref{thm:Bound_EGW1} and given in Appendix \ref{appen:proof_Bound_EGW2}.
\begin{theorem}[EOT cost neural estimation; bound 2] \label{thm:Bound_EGW2}
    Let $\epsilon>0$ and set $m_k=\log k\vee 1$. Assuming $k$ is sufficiently large, we have
\begin{equation}
\label{Bound_EGW2}
\begin{aligned}
   &\sup _{(\mu, \nu) \in \mathcal{P}(\cX) \times \mathcal{P}(\cY)} \mspace{-5mu}\mathbb{E}\left[\left|\widehat\OT_{k,m_k}^{\mspace{1mu}\eps}(X^n, Y^n)\mspace{-2mu}-\mspace{-2mu}\OT^{\eps}(\mu, \nu)\right|\right]\\
   &\quad\quad\quad\quad\quad\quad\lesssim_{d} \left(1+\frac{1}{\eps^{\left\lfloor\frac{d}{2}\right\rfloor+2}}\right) k^{-\frac{1}{2}}+\min\left\{\left(1+\frac{1}{\eps^{\left[\frac{d}{2}\right\rceil}}\right)(\log k)^2\,, \sqrt{k}\log k\right\} n^{-\frac{1}{2}}.
\end{aligned}
\end{equation}
\end{theorem}
\begin{remark}[NN size]
    We can provide a partial account of the requirement that $k$ is large enough. Specifically, for the bound to hold we need $k$ to satisfy $\log k \geq C(1+\eps^{1-s})$, where $C$ is the constant from Theorem \ref{thm:Bound_EGW1}. It is, however, challenging to quantify the exact threshold on $k$ required for the theorem to hold due to the implicit nature of $C$.
\end{remark}

Lastly, we move to account for the quality of the neural plan that is induced by the EOT NE (see \eqref{eq:NE_plan}) by comparing it, in KL divergence, to the true EOT plan $\pi^\eps_\star$.
\begin{theorem}[EOT alignment plan neural estimation]
\label{thm:Bound_NeuralCoupling}
    Suppose that $\mu \in \mathcal{P}_{\mathsf{ac}}(\cX)$. Let $\hat f_\star$ be a maximizer of $\widehat\OT^\eps_{k,a}(X^n,Y^n)$ from \eqref{eq:NE_EOT}, with $a$ as defined in Theorem \ref{thm:Bound_EGW1}. Then, the induced neural plan $\pi^\eps_{\hat f_\star}$ from \eqref{eq:NE_plan} satisfies
    \begin{equation}
    \label{Bound_NeuralCoupling}
    \begin{aligned}
&\mathbb{E}\left[\KL\left(\pi_\star^\eps\middle\|\pi^\eps_{\hat f_\star}\right)\right]\\
&\quad\quad\quad\quad\quad\quad\lesssim_{d}\left(1+\frac{1}{\eps^{\left\lfloor \frac{d}{2}\right\rfloor+3}}\right) k^{-\frac{1}{2}}+\min \left\{1+\frac{1}{\eps^{\left\lceil \frac{3d}{2}\right\rceil+5}},\left(1+\frac{1}{\eps^{\left\lfloor\frac{d}{2}\right\rfloor+3}}\right) \sqrt{k}\right\}n^{-\frac{1}{2}} .
    \end{aligned}
    \end{equation}
where $\pi^\eps_*$ is optimal coupling of EOT problem \eqref{EOT}.
\end{theorem}

Theorem \ref{thm:Bound_NeuralCoupling} is proved in Appendix \ref{appen:proof_thm3}. The key step in the derivation shows that the KL divergence between the alignment plans, in fact, equals the gap between the EOT cost $\OT^\eps$ and its neural estimate from \eqref{eq:NE_EOT}, up to a multiplicative $\eps$ factor. Having that, the result follows by invoking Theorem~\ref{thm:Bound_EGW1}.

\begin{remark}[Extension to Sigmoidal NNs]
\label{rem:sigmoid_bound}
The results of this section readily extend to cover sigmoidal NNs, with a slight modification of some parameters. Specifically, one has to replace $s$ from Theorem \ref{thm:Bound_EGW1} with $\tilde s=\lfloor d / 2\rfloor+2$ and consider the sigmoidal NN class, with nonlinearity $\psi(z)=\left(1+e^{-z}\right)^{-1}$ (instead of ReLU) and parameters satisfying
\begin{align*}
\max _{1 \leq i \leq k}\left\|w_i\right\|_1 \vee\left|b_i\right| \leq k^{\frac{1}{2}}\log k,\ \max _{1 \leq i \leq k}\left|\beta_i\right| \leq 2ak^{-1},\ \left|b_0\right| \leq a,\ \left\|w_0\right\|_1 = 0.
\end{align*}
The proofs of Theorems \ref{thm:Bound_EGW1}-\ref{thm:Bound_NeuralCoupling} then go through using the second part of Proposition 10 from \cite{sreekumar2022neural}, which relies on controlling the so-called Barron coefficient (cf. \cite{barron1992neural,barron1993universal,yukich1995sup}).
\end{remark}

\begin{remark}[Convergence rates of Sinkhorn's algorithm]
Neural estimation is proposed as a more scalable alternative to Sinkhorn's algorithm for computing the EOT cost/plan, e.g., by enabling the usage of mini-batches. We comment here on the rate of convergence that the Sinkhorn-based approach achieves. Denote the output of the Sinkhorn algorithm running on empirical measures, each over $n$ samples, by $\widetilde\OT^\eps(X^n,Y^n)$. The effective error can be decomposed as:
\begin{align*}
|\widetilde\OT^\eps(X^n,Y^n)-\OT(\mu,\nu)|\leq |\OT(\mu,\nu)-\OT(\mu_n,\nu_n)|+|\OT(\mu_n,\nu_n)-\widetilde\OT^\eps(\mu_n,\nu_n)|,
\end{align*}
where first term decays as $O(n^{-\frac{1}{2}})$ \cite{mena2019statistical}, while the second exhibits a convergence rate of $o_P(n^{-1})$ within $o_P(\log (n \log (n)))$ iterations \cite{goldfeld2024limit}.
\end{remark}

\section{Numerical Experiments}
This section illustrates the performance of the EOT distance neural estimator via experiments with synthetic data. Specifically, we compute the estimate from \eqref{eq:NE_EOT} under various settings, allowing $a$ to be  unrestricted so as to enable optimization over the whole parameter space. We train the parameters of the ReLU network using the Adam algorithm \cite{kingma2014adam}. We use an epoch number of 20, learning rate $10^{-3}$ and choose a best batch size from $\{2,4,8,16,32,64,128\}$. We test our EOT distance neural estimator by estimating the EOT cost and optimal plan between uniform and Gaussian distribution in different dimensions. We consider dimensions $d\in \{1,16,64,128\}$, and for each $d$, employ a ReLU network of size $k\in \{16,64,128,256\}$, respectively. Accuracy is measured using the relative error $\big|\widehat{\OT}_{k, a}^{\eps}\left(X^n, Y^n\right)-\widetilde\OT^\eps(\mu,\nu)\big|/\,\widetilde\OT^\eps(\mu,\nu)$, where $\widetilde\OT^\eps(\mu,\nu)$ is regarded as the ground truth, which we obtain by running Sinkhorn algorithm \cite{cuturi2013sinkhorn} with $n=10,000$ samples (which we treat as $n\to\infty$ as it is $\times 5$ more than the largest sample set we use for our neural estimator). Each of the presented plots is averaged over 20 runs.

\begin{figure}[!t]
  \begin{subfigure}{0.345\textwidth}
    \centering
    \includegraphics[width=\textwidth]{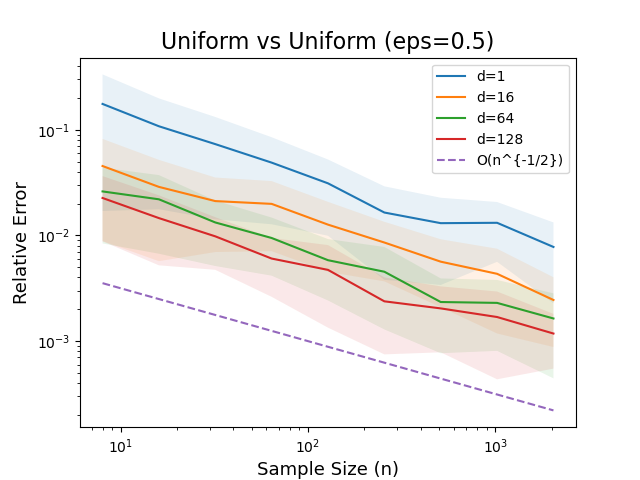}
    \caption{}
    \label{fig:uniform}
  \end{subfigure}%
  \hspace{-6mm}
  \begin{subfigure}{0.345\textwidth}
    \centering
    \includegraphics[width=\textwidth]{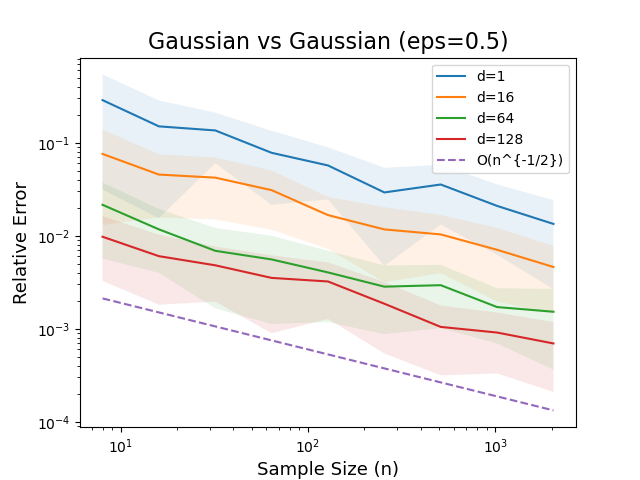}
    \caption{}
    \label{fig:gaussian}
  \end{subfigure}
  \hspace{-6mm}
    \begin{subfigure}{0.345\textwidth}
    \centering
    \includegraphics[width=\textwidth]{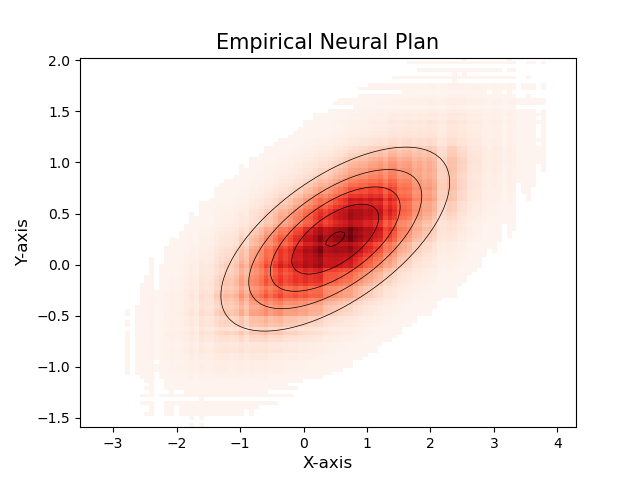}
    \caption{}
    \label{fig:Neural_Plan}
  \end{subfigure}%
  \caption{Neural Estimation of EOT distance: (a) Relative error for the case where $\mu=\nu=\mathrm{Unif}\big([-1/\sqrt{d},1/\sqrt{d}]^d\big)$; (b) Relative error for $\mu,\nu$ as Gaussian distributions with randomly generated mean vectors and covariance matrices; (c) Learned neural plan (in red) versus the true optimal EOT optimal plan (whose density is represented by the back contour lines).}
  \label{fig:simulated_data}
\end{figure}

We first consider the EOT distance with $\eps=0.5$ between two uniform distribution over a hypercube, namely,  $\mu=\nu=\mathrm{Unif}\big([-1/\sqrt{d},1/\sqrt{d}]^d\big)$. Figure \ref{fig:uniform} plots the EOT neural estimation error versus the sample size $n\in\{8,16,32,64,128,256,512,1024,2048\}$ in a log-log scale. The curves exhibit a slope of approximately $-1/2$ for all dimensions, which validates our theory. Notably, this rate is uniform across dimensions, like the bounds from Theorems \ref{thm:Bound_EGW1} and \ref{thm:Bound_EGW2} suggest.

Next, we test the EOT NE on unbounded measures. To that end, we set $\eps = 0.5$ and take $\mu,\nu$ as $d$-dimensional Gaussian distributions with randomly generated mean vectors and covariance matrices. Specifically, the mean vectors are randomly sampled from $d$-dimensional standard Gaussian, while the two covariance matrices are of the form $\mathbf{B}^{\intercal}\mathbf{B}+1/(3d)\mathbf{I}_d$, where $\mathbf{I}_d$ is a $d\times d$ identity and $\mathbf{B}$ is a matrix whose entries are randomly sampled from $\mathrm{Unif}([-1/d,1/d])$. Note that the generated covariance matrix is positive semi-definite with eigenvalues set to lie in $[\frac{1}{3d},\frac{1}{d}]$. Figure \ref{fig:gaussian} plots the relative EOT neural estimation error for this Gaussian setting, again showing a parametric convergence rate for all considered dimensions. 

Lastly, we assess the quality of the neural plan learned from our NE. Since doing so requires knowledge of the true (population) optimal plan $\pi^\eps_\star$, we consider the EOT distance between Gaussians, for which a closed form expression for the optimal plan was derived in \cite{janati2020entropic}. We take $\eps=0.5$, $\mu=\mathcal{N}\left(0.5,1\right)$, and $\nu=\mathcal{N}(0.25,0.25)$. By Theorem 3.1 of \cite{le2022entropic}, the optimal EOT plan is given by
$$
\pi^\eps_{\star} \sim \mathcal{N}\left(\left(\begin{array}{l}
0.5 \\
0.25
\end{array}\right),\left(\begin{array}{cc}
1 & \frac{\sqrt{5}-1}{4} \\
\frac{\sqrt{5}-1}{4} & 0.25
\end{array}\right)\right)
$$
Figure \ref{fig:Neural_Plan} compares the neural coupling learned from our algorithm, shown in red, to the optimal $\pi^\eps_\star$ given above, whose density is represented by the black contour lines. The neural coupling is learned using $n=10^4$ samples and is realized by a NN with $k=32$ neurons. There is a clear correspondence between the two, which supports the result of Theorem \ref{thm:Bound_NeuralCoupling}.

\section{Conclusion}

This work proposed a novel neural estimation technique for the EOT distance with quadratic costs between Euclidean mm spaces. The estimator leveraged the semi-dual formulation of EOT. Our approach yielded estimates not only for the EOT distance value but also for the optimal plan. Non-asymptotic formal guarantees on the quality of the NE were provided, under the sole assumption of compactly supported population distributions, with no further regularity conditions imposed. Our bounds revealed optimal scaling laws for the NN and the dataset sizes that ensure parametric (and hence minimax-rate optimal) convergence. The proposed estimator was tested via numerical experiments on synthetic data, demonstrating its accuracy, scalability, and fast convergence rates that match the derived theory. 

Future research directions stemming from this work are abundant. First, our theory currently accounts for NEs realized by shallow NNs, but deep nets are oftentimes preferable in practice. Extending our results to deep NNs should be possible by utilizing existing function approximation error bounds \cite{bresler2020sharp}, although these bounds may not be sharp enough to yield the parametric rate of convergence. Another limitation of our analysis is that it requires compactly supported distributions. It is possible to extend our results to distributions with unbounded supports using the technique from \cite{sreekumar2022neural} that considers a sequence of restrictions to balls of increasing radii. Unfortunately, as in \cite{sreekumar2022neural}, rate bounds obtained from this technique would be sub-optimal. Obtaining sharp rates for the unboundedly supported case would require new ideas and forms an interesting research direction. Lastly, while EOT serves as an important approximation of OT, neural estimation of the OT distance itself is a challenging and appealing research avenue. One may attempt to directly approximate this objective by NNs, but dual OT potential generally lack sufficient regularity to allow quantitative approximation bounds. Assuming smoothness of the population distributions, and employing estimators that adapt to this smoothness, e.g., based on kernel density estimators or wavelets \cite{deb2021rates,manole2021plugin}, may enable deriving sharp rates of convergence. 

\bibliographystyle{alpha}
\bibliography{ref}

\newpage
\appendix
\section{Proofs}
\label{appen:Proofs}

We first introduce a technical result from approximation theory that will be used in the subsequent derivations. The following result, which is a restatement of Proposition 10 from \cite{sreekumar2022neural}, states that a sufficiently smooth function over a compact domain can be approximated to within $O(k^{-1/2})$ error by a shallow NN.
\begin{proposition}[Approximation of smooth functions; Proposition 10 from \cite{sreekumar2022neural}]
\label{prop:approx_restated}
Let $\cX\subseteq\RR^d$ be compact and $g:\cX \rightarrow \RR$. Suppose that there exists an open set $\cU \supset  \cX$, $b \geq 0$, and $ \tilde g \in  \cC_b^{s_{\mathsf{KB}}}(\cU)$, $s_{\mathsf{KB}}\coloneqq \lfloor d/2\rfloor+3$, such that $g= \tilde g\big|_\cX$.  
Then, there exists $f \in \cF_{k,d}\left(\bar c_{b,d,\|\cX\|}\right)$, where $\bar c_{b,d,\|\cX\|}$ is given in Equation (A.15) of \cite{sreekumar2022neural}, such~that 
\[\|f-g\|_\infty\lesssim \bar{c}_{b,d,\|\cX\|}d^{\frac 12} k^{-\frac 12}.\] 
\end{proposition}  

This proposition will allow us to control the approximation error of the EOT NE. To invoke
it, we will establish smoothness of the semi-dual EOT potentials (see Lemma \ref{lemma:regularity_EOT} ahead). The smoothness of potentials stems from the
presence of the entropic penalty and the smoothness of the quadratic cost function.

\subsection{Proof of Theorem \ref{thm:Bound_EGW1}}
\label{appen:proof_Bound_EGW1}

 For $(\mu, \nu) \in \mathcal{P}(\cX) \times \mathcal{P}(\cY)$, define the population-level neural EOT cost as
\begin{equation}
\OT^{\eps}_{k,a}(\mu, \nu)\coloneqq\sup_{f\in \cF_{k,d}(a)}\int f d\mu+\int f^{c,\eps}d\nu.\label{eq:neural_pop_EOT}
\end{equation}
We decompose the neural estimation error into the approximation and empirical estimation errors:
\begin{equation}
    \begin{aligned}
\label{eq:EGW_NE_error_decomposition}
&\mathbb{E}\left[\left|\widehat\OT_{k,a}^{\mspace{1mu}\eps}\left(X^n, Y^n\right)-\OT^{\eps}(\mu, \nu)\right|\right]\nonumber\\
&\qquad\qquad\leq\underbrace{\left|\OT_{k,a}^{\mspace{1mu}\eps}(\mu, \nu)-\OT^{\eps}(\mu, \nu)\right|}_{\text{ Approximation error }}\nonumber+\underbrace{\mathbb{E}\left[\left|\OT^{\eps}_{k,a}(\mu, \nu)-\widehat{\OT}^{\eps}_{k,a}\left(X^n, Y^n\right)\right|\right]}_{\text{ Estimation error }},
\end{aligned}
\end{equation}

and analyze each term separately. 

\medskip

\noindent\underline{Approximation error.} Proposition \ref{prop:approx_restated} provides a sup-norm approximation error bound of a smooth function by a NN. To invoke it, we first study the semi-dual EOT potentials and show that they are indeed smooth functions, i.e., admit an extension to an open set with sufficiently many bounded derivatives. The following lemma establishes regularity of semi-dual potentials for $\OT^{\eps}(\mu,\nu)$; after stating it we shall account for the extension. 

\begin{lemma}[Uniform regularity of EOT potentials]\label{lemma:regularity_EOT}
There exist semi-dual EOT potentials $(\varphi,\varphi^{c,\eps})$ for $\OT^{\mspace{1mu}\eps}(\mu, \nu)$, such that
\begin{equation}
\begin{aligned}
\label{eq:regularity_EOT}
    \|\varphi\|_{\infty,\cX} &\leq 2d\\
    \left\|D^\alpha \varphi\right\|_{\infty,\cX} &\leq C_{s}\left(1+\varepsilon^{1-s}\right)\left(1+2 \sqrt{d}\right)^{s} \text { with } 1 \leq|\alpha| \leq s,
\end{aligned}
\end{equation}
for any $s \geq 2$ and some constant $C_s$ that depends only on $s$. Analogous bounds hold for $\varphi^{c,\eps}$.
\end{lemma}
The lemma is proven in Appendix \ref{appen:proof:lemma:regularity_EOT}. The derivation is similar to that of Lemma 4 in \cite{zhang2022gromov}, but the bounds are adapted to the compactly supported case and present an explicit dependence on $\eps$ (as opposed to the $\eps=1$ assumption that was imposed in that work). 

Let $(\varphi,\varphi^{c,\eps})$ be semi-dual potentials as in  Lemma \ref{lemma:regularity_EOT} (i.e., satisfying \eqref{eq:regularity_EOT}) with the normalization $\int \varphi d \mu=\int \varphi^{c,\eps} d \nu=\frac{1}{2} \mathrm{OT}^{\eps}(\mu, \nu)$. Define the natural extension of $\varphi$ to the open ball of radius $\sqrt{2d}$:
\[
\tilde\varphi(x)\coloneqq \mspace{-1mu}-\eps \log\mspace{-4mu} \int \mspace{-4mu} \exp\mspace{-1mu} \left(\mspace{-1mu}\frac{\psi(y)-c(x, y)}{\eps}\mspace{-1mu}\right)\mspace{-1mu} d \nu(y), x \in B_d(\mspace{-1mu}\sqrt{2d}),
\]
and notice that $\tilde\varphi\big|_\cX=\varphi$, pointwise on $\cX$. Similarly, consider its $(c,\eps)$-transform $\tilde\varphi^{c,\eps}$ extended to $ B_d(\sqrt{2d})$, and again observe that $\tilde\varphi^{c,\eps}\big|_\cY=\varphi^{c,\eps}$. Following the proof of Lemma \ref{lemma:regularity_EOT}, one readily verifies that for any $s\geq 2$, we have
\begin{equation}
    \begin{aligned}
    \label{derivative_est_of_extension}
\|\tilde{\varphi}\|_{\infty,B_{d}\left(\sqrt{2d}\right)} & \mspace{-1mu}\leq \mspace{-1mu}4d \\
\left\|D^\alpha \tilde{\varphi}\right\|_{\infty,B_{d}\left(\sqrt{2d}\right)} & \mspace{-1mu}\leq\mspace{-1mu} C_s\left(1+\varepsilon^{1-s}\right)(1+2 \sqrt{2d})^s, 1 \leq|\alpha| \leq s.
\end{aligned}
\end{equation}
Recall that $s_{\mathsf{KB}}=\left\lfloor d / 2\right\rfloor+3$ and set
\begin{equation}
    \label{eq:constant_C_dx_dy}
    C_d\coloneqq (1+2 \sqrt{2d})^{s_{\mathsf{KB}}}.
\end{equation}
By \eqref{derivative_est_of_extension} and \eqref{eq:constant_C_dx_dy}, we now have
\begin{equation}
\label{eq:constant_b}
  \max _{\alpha:\,|\alpha| \leq s_{\mathsf{KB}}}\left\|D^\alpha \tilde{\varphi}\right\|_{\infty, B_{d}\left(\sqrt{2d}\right)} \leq C_{s_{\mathsf{KB}}}C_d(1+\eps^{1-s_{\mathsf{KB}}}) \coloneqq b,
\end{equation}
and so $\tilde{\varphi} \in \cC_b^{s_{\mathsf{KB}}}\big(B_d(\sqrt{2d})\big)$.

Noting that $\cX\subset B_{d}(\sqrt{2d})$, by Proposition \ref{prop:approx_restated}, there exists $f \in \cF_{k,d}\left(\bar c_{b,d}\right)$ such that 
\begin{equation}
\label{eq:EOTpotential_approx}
  \left\|\varphi-f\right\|_{\infty,\cX}\lesssim \bar{c}_{b,d}d^{\frac 12} k^{-\frac 12},  
\end{equation}
where $\bar{c}_{b,d} = b\, \bar{c}_{d}$ and $\bar{c}_{d}$ is defined as (see \cite[Equation (A.15)]{sreekumar2022neural})
\begin{equation}
\label{eq:constant_c_bar}
    \begin{aligned}
\bar{c}_{d} & \coloneqq \left(\kappa_{d}{d}^{\frac{3}{2}} \vee 1\right)\pi^{\frac{d}{2}} \Gamma\left(\frac{d}{2}+1\right)^{-1}(\operatorname{rad}(\cX)+1)^{d} \\
& \quad\quad\quad\times  2^{s_{\mathsf{KB}}} d\left(\frac{1-{d}^{\frac{s_{\mathsf{KB}}}{2}}}{1-\sqrt{d}}+{d}^{\frac{s_{\mathsf{KB}}}{2}}\right) \max _{\|\alpha\|_1 \leq s_{\mathsf{KB}}}\left\|D^\alpha \Psi\right\|_{\infty, B_d(0.5)},
\end{aligned}
\end{equation}
with $\kappa_{d}^2\coloneqq \left(d+{d}^{(s_{\mathsf{KB}}-1)}\right) \int_{\mathbb{R}^d}\left(1+\|\omega\|^{2\left(s_{\mathsf{KB}}-2\right)}\right)^{-1} d \omega$, $\operatorname{rad}(\cX)=0.5 \sup _{x, x^{\prime} \in \cX}\left\|x-x^{\prime}\right\|$ and $\Psi(x) \propto \exp \left(-\frac{1}{0.5-\|x\|^2}\right) \mathbbm{1}_{\{\|x\|<0.5\}}$ as the canonical mollifier normalized to have unit mass. 

\medskip
Our last step is to lift the sup-norm neural approximation bound on the semi-dual potential from \eqref{eq:EOTpotential_approx} to a bound on the approximation error of the corresponding EOT cost. The following lemma is proven in Appendix \ref{appen:proof_neural_approx_reduction}.
\begin{lemma}[Neural approximation error reduction]
\label{lemma:neural_approx_reduction}
Fix $(\mu,\nu) \in \mathcal{P}(\cX) \times \mathcal{P}(\cY)$ and let $\varphi$ be the semi-dual EOT potential for $\OT^\eps(\mu,\nu)$ from Lemma \ref{lemma:regularity_EOT}. For any $f \in \cF_{k,d}\left(a\right)$, we have
\[
    \left|\OT^{\eps}(\mu,\nu)-\OT^{\eps}_{k,a}(\mu,\nu)\right| \leq  2\left\|\varphi-f\right\|_{\infty, \cX}.
\]
\end{lemma}

Setting $a=\bar{c}_{b,d}$ and combining Lemma \ref{lemma:neural_approx_reduction} 
 with \eqref{eq:EOTpotential_approx}, we obtain
\begin{equation}
    \label{eq:EGW_approx_error}
    \left|\OT^{\eps}_{k,a}(\mu,\nu)-\OT^\eps(\mu,\nu)\right|\lesssim_{d}\left(1+\frac{1}{\eps^{\left\lfloor \frac{d}{2}\right\rfloor+2}}\right) k^{-\frac{1}{2}}.
\end{equation}

\medskip

\noindent\underline{Estimation error}. Set $\cF^{c,\eps}(a)\coloneqq \left\{f^{c,\eps}: f \in \cF_{k,d}(a)\right\}$, and first bound
\begin{align}\label{eq:emp_error_decomposition}
    &\mathbb{E}\left[ \left|\OT^\eps_{k,a}(\mu, \nu)-\widehat\OT^\eps_{k,a}\left(X^n, Y^n\right)\right|\right]\nonumber\\
    & \leq\underbrace{n^{-\frac{1}{2}} \mathbb{E}\left[\sup_{f \in \cF_{k,d}(a)} \mspace{-3mu}n^{-\frac{1}{2}} \left| \sum_{i=1}^n \big(f(X_i)-\mathbb{E}_{\mu}[f]\big) \right|  \right]}_{(\mathrm{I})}+\underbrace{n^{-\frac{1}{2}}\mathbb{E}\left[\sup_{f \in \cF^{c,\eps}(a)}\mspace{-3mu}n^{-\frac{1}{2}}\left| \sum_{j=1}^n  \big(f(Y_j)-\mathbb{E}_{\nu}[f]\big) \right|  \right]}_{(\mathrm{II})}.
\end{align}
To control these expected suprema, we again require regularity of the involved function, as stated in the next lemma.

\begin{lemma}
\label{lemma:c_transform_regularity}
    Fix $c \in \cC^{\infty}$, the $(c,\eps)$-transform of NNs class $\cF_{k, d}(a)$ satisfies the following uniform smoothness properties:
\[
\max _{\alpha:|\alpha|_1 \leq s}\left\|D^\alpha f^{c,\eps}\right\|_{\infty, \cY} \leq R_s(1+a)(1+\eps^{1-s})
\]
for any $s \geq 2$, $f \in \cF_{k, d}(a)$ and some constant $R_s$ that depends only on $s,d$.
\end{lemma}
The only difference between Lemmas \ref{lemma:regularity_EOT} and \ref{lemma:c_transform_regularity} is that here we consider the $(c,\eps)$-transform of NNs, rather than of dual EOT potentials. As our NNs are also compactly supported and bounded, the derivation of this result is all but identical to the proof of Lemma \ref{lemma:regularity_EOT}, and is therefore omitted to avoid repetition.

We proceed to bound Terms $(\mathrm{I})$ and $(\mathrm{II})$ from \eqref{eq:emp_error_decomposition}. For the first, consider
\begin{align}
    \mathbb{E}\left[\sup_{f \in \cF_{k,d}(a)} n^{-\frac{1}{2}} \left| \sum_{i=1}^n \left(f(X_i)-\mathbb{E}_{\mu}[f]\right) \right|  \right]
   & \stackrel{(a)}{\lesssim} \mathbb{E}\left[\int_0^{\infty} \sqrt{\log N\left(\delta, \cF_{k,d}(a),\|\cdot\|_{2, \mu_n}\right)} d \delta\right]\nonumber\\
    & \leq \int_0^{\infty} \sqrt{\sup _{\gamma \in \mathcal{P}(\cX)} \log N\left(\delta, \cF_{k,d}(a),\|\cdot\|_{2, \gamma}\right)} d \delta\nonumber\\
    & \stackrel{(b)}{=} \int_0^{12a} \sqrt{\sup _{\gamma \in \mathcal{P}(\cX)} \log N\left(\delta, \cF_{k,d}(a),\|\cdot\|_{2, \gamma}\right)} d \delta\nonumber\\
    & \lesssim a  \int_0^1 \sqrt{\sup _{\gamma \in \mathcal{P}(\cX)} \log N\left(6a \delta, \cF_{k,d}(a),\|\cdot\|_{2, \gamma}\right)}d \delta\nonumber\\
    & \stackrel{(c)}{\lesssim}ad^{\frac{3}{2}},\label{eq:term1_bound}
\end{align}
where:

(a) follows by \cite[Corollary 2.2.8]{van1996springer} since $n^{-\frac{1}{2}} \sum_{i=1}^n \sigma_i f\left(X_i\right)$, where $\{\sigma_i\}_{i=1}^n$ are i.i.d Rademacher random variables, is sub-Gaussian w.r.t. pseudo-metric $\|\cdot\|_{2, \mu_n}$ (by Hoeffding's inequality) ;

(b) is since $\bar{C}\left(\left|\cF_{k,d}(a)\right|, \cX\right) \leq 3 a(\|\cX\|+1)=6a$ and $N\left(\delta, \cF_{k,d}(a),\|\cdot\|_{2, \gamma}\right)=1$, whenever $\delta > 12a$;

(c) uses the bound $$\int_0^1 \sqrt{\sup _{\gamma \in \mathcal{P}(\cX)} \log N\left(6a \delta, \cF_{k,d}(a),\|\cdot\|_{2, \gamma}\right)} d \delta \lesssim d^{\frac{3}{2}},$$
which follows from step (A.33) in \cite{sreekumar2022neural}.

\medskip
For Term $(\mathrm{II})$, let $s=\left\lceil d  / 2\right\rceil+1$, and consider
\begin{align}
    \mathbb{E}\left[\sup _{f \in \cF^{c,\eps}(a)} n^{-\frac{1}{2}}\left|\sum_{i=1}^n\big(f\left(Y_i\right)-\mathbb{E}_\mu[f]\big)\right|\right]
    &\stackrel{(a)}{\lesssim} \int_0^{12 a} \sqrt{\sup _{\gamma \in \mathcal{P}(\cY)} \log N\left(\delta, \cF^{c,\eps}(a),\|\cdot\|_{2, \gamma}\right)} d \delta\nonumber\\
    &\lesssim \int_0^{12 a} \sqrt{\sup _{\gamma \in \mathcal{P}(\cY)} \log N_{[\,]}\left(2\delta, \cF^{c,\eps}(a),\|\cdot\|_{2, \gamma}\right)} d \delta\nonumber\\
    & \stackrel{(b)}{\lesssim} K_s \int_{0}^{12 a} \left(\frac{R_s(1+a)(1+\eps^{1-s})}{2\delta}\right)^{\frac{d}{2s}}d\delta\nonumber\\
    & \lesssim_{d} a(1+a)(1+\eps^{1-s}),\label{eq:term2_bound1}
\end{align}
where $R_s$ is the constant from Lemma \ref{lemma:c_transform_regularity} (which depends only $s,d$), (a) follows by a similar argument to that from the bound on Term $(\mathrm{I})$, along with equation \eqref{eq:lipshcitz_c_tranform}, which specifies the upper limit for entropy integral, while (b) follows by Lemma \ref{lemma:c_transform_regularity} and \cite[Corollary 2.7.2]{van1996springer}, which upper bounds the bracketing entropy number of smooth functions on a bounded convex support. 

To arrive at the effective error bound from Theorem \ref{thm:Bound_EGW1}, we provide a second bound on Term $(\mathrm{II})$. This second bound yields a better dependence on dimension (namely, only the smaller dimension $d$ appears in the exponent) at the price of another $\sqrt{k}$ factor. Neither bound is uniformly superior over the other, and hence our final result will simply take the minimum of the two. By \eqref{eq:lipshcitz_c_tranform} from the proof of Lemma \ref{lemma:neural_approx_reduction}, we have 
\[
N\left(\delta, \mathcal{F}^{c,\eps}(a),\|\cdot\|_{2, \gamma}\right) \leq N\left(\delta, \cF_{k, d}(a),\|\cdot\|_{\infty,\cX}\right).
\]
Invoking Lemma 2 from \cite{sreekumar2022neural}, which upper bounds the metric entropy of ReLU NNs class on the RHS above, we further obtain
\begin{equation}
\begin{aligned}
 &\log N\left(\delta, \cF_{k, d}(a),\|\cdot\|_{\infty, \cX}\right) \\
 &\leq \big((d+2) k+d+1\big)\log\left(1+20 a(\|\cX\|+1) \delta^{-1}\right),
\end{aligned}
\end{equation}
and proceed to bound Term $(\mathrm{II})$ as follows:
\begin{align}
     \mathbb{E}\left[\sup_{f \in \mathcal{F}^{c,\eps}}n^{-\frac{1}{2}}\left| \sum_{j=1}^n  \big(f(Y_j)-\mathbb{E}_{\nu}[f]\big) \right|  \right] 
     &\lesssim a\int_0^1 \sqrt{\sup _{\gamma \in \mathcal{P}(\cX)} \log N\left(6 a \delta, \mathcal{F}^{c,\eps}(a),\|\cdot\|_{2, \gamma}\right)} d \delta  \nonumber \\
     & \lesssim ad^{\frac{1}{2}}\sqrt{k}\int_0^1 \sqrt{\log(1+7\delta^{-1})} d \delta\nonumber\\
     & \lesssim ad^{\frac{1}{2}}\sqrt{k}.\label{eq:term2_bound2}
\end{align}

\medskip

Inserting \eqref{eq:term1_bound}, \eqref{eq:term2_bound1}, and \eqref{eq:term2_bound2} back into \eqref{eq:emp_error_decomposition}, we obtain the desired bound on the empirical estimation error by setting $a=\bar{c}_{b,d}$, as was defined in approximation error analysis:
\begin{equation}
\begin{aligned}
    \label{eq:EGW_est_error}
    &\mathbb{E}\left[\left|\OT^\eps_{k,a}(\mu, \nu)-\widehat\OT^\eps_{k,a}\left(X^n, Y^n\right)\right|\right]\lesssim_{d} \min\left\{1+\frac{1}{\eps^{\left\lceil \frac{3d}{2}\right\rceil+4}}\,,\left(1+\frac{1}{\eps^{\left\lfloor \frac{d}{2}\right\rfloor+2}}\right)\sqrt{k}\right\} n^{-\frac{1}{2}}.
\end{aligned}
\end{equation}

The proof is concluded by plugging the approximation error bound from \eqref{eq:EGW_approx_error} and the estimation error bounds from \eqref{eq:EGW_est_error} into \eqref{eq:EGW_NE_error_decomposition}, and supremizing over $(\mu, \nu) \in \mathcal{P}(\cX) \times \mathcal{P}(\cY)$, while noting that all the above bounds holds uniformly in the two distributions. 
\qed

\subsection{Proof of Theorem \ref{thm:Bound_EGW2}}
\label{appen:proof_Bound_EGW2}

This proof is similar to that of Theorem \ref{thm:Bound_EGW1}, up to minor modifications. For brevity, we only highlight the required changes. Note that for $k$ with $m_k \geq \bar{c}_{b,d}$, where the latter is defined in proof of Theorem \ref{thm:Bound_EGW1} (see \eqref{eq:constant_b} and \eqref{eq:constant_c_bar}), we have $\cF_{k,d}(\bar{c}_{b,d})\subset\cF_{k,d}(m_k)$. Hence, by Lemma \ref{lemma:neural_approx_reduction} and \eqref{eq:EOTpotential_approx}, there exists a NN $f\in\cF_{k,d}(\bar{c}_{b,d})$, such that
\begin{align*}
   \left|\OT^{\eps}(\mu, \nu)-\OT^{\eps}_{k,m_k}(\mu, \nu)\right|  &\leq 2\left\|\varphi-f\right\|_{\infty,\cX}\lesssim_{d} \left(1+\frac{1}{\eps^{\left\lfloor \frac{d}{2}\right\rfloor+2}}\right) k^{-\frac{1}{2}}. 
\end{align*}

Next, for estimation error, by setting $a=m_k\coloneqq\log k\vee 1$ in  \eqref{eq:term1_bound}, \eqref{eq:term2_bound1}, and \eqref{eq:term2_bound2} (instead of $a=\bar{c}_{b,d}$ as in the proof of Theorem \ref{thm:Bound_EGW1}), we arrive at
\begin{align*}
&\mathbb{E}\left[ \left|\OT^\eps_{k,m_k}(\mu, \nu)-\widehat\OT^\eps_{ k,m_k}\left(X^n, Y^n\right)\right|\right] \lesssim_{d}\min\left\{\left(1+\frac{1}{\eps^{\left[\frac{d}{2}\right\rceil}}\right)(\log k)^2\,, \sqrt{k}\log k\right\} n^{-\frac{1}{2}}.
\end{align*}
Combining both bounds completes the proof.
\qed

\subsection{Proof of Theorem \ref{thm:Bound_NeuralCoupling}}
\label{appen:proof_thm3}
Define $\Gamma(f): = \int_{\cX} fd\mu+\int_{\cY} f^{c, \eps}d \nu$, and let $\varphi_\star$ be optimal potential of $\OT^{\eps}(\mu,\nu)$, solving semi-dual formulation. Denote the corresponding optimal coupling by $\pi^\eps_\star$. We first show that for any continuous $f:\cX\to\RR$, the following holds:
\begin{equation}
\label{eq:stability_coupling}
    \Gamma\left(\varphi_\star
    \right)-\Gamma\left(f\right) =\eps\KL\big(\pi^\eps_\star\big\|\pi^\eps_{f}\big),
\end{equation}
where (see \eqref{eq:NE_plan}) 
\[
d \pi^\eps_{f}(x, y) = \frac{\exp \left(\frac{f(x)-c(x, y)}{\eps}\right)}{\int_{\cX} \exp \left(\frac{f-c(\cdot, y)}{\eps}\right) d\mu} d\mu\otimes\nu(x,y).
\]

The derivation is inspired by the proof of \cite[Theorem 2]{mokrov2023energy}, with several technical modifications. 
Since $\mu \in \mathcal{P}_{\mathsf{ac}}(\cX)$ with Lebesgue density $\frac{d\mu}{dx}$, define its energy function $E_{\mu}: \cX \rightarrow \mathbb{R}$ by $\frac{d \mu(x)}{d x} \propto \exp \left(-E_{\mu}(x)\right)$. Also define conditional distribution $d\pi^\eps_{f}(\cdot| y) \coloneqq  \frac{d\pi_{f}(\cdot,y)}{d\nu(y)}$, and set $\tilde{f}\coloneqq f-\eps E_{\mu}(x)$. We have
\begin{align*}
    \frac{d\pi^\eps_{f}(x| y) }{d x}& =\frac{d\pi^\eps_{f}(x,y)}{d\nu(y)dx}\mspace{-.5mu}=\mspace{-.5mu} \frac{\exp \left(\frac{f(x)-c(x, y)}{\eps}\right)\frac{d\mu(x)}{dx}}{\int_{\cX} \exp \left(\frac{f(x')-c(x', y)}{\eps}\right)\frac{d\mu}{dx'}(x')dx'}\mspace{-.5mu}=\mspace{-.5mu} \frac{\exp \left(\frac{\tilde{f}(x)-c(x, y)}{\eps}\right)}{\int_{\cX} \exp \left(\frac{\tilde{f}(x')-c(x', y)}{\eps}\right)dx'}.
\end{align*}  
Define the shorthands $F_{f}(y)\coloneqq\int_{\cX} \exp \left(\frac{f(x)-c(x, y)}{\eps}\right)dx$ and $Z \coloneqq \int_\cX \exp \big(-E_{\mu}(x)\big)dx$, and note that the $(c,\eps)$-transform of $f$ can be expressed as
\begin{align*}
   f^{c, \eps}(y)&=-\eps \log \left(\int_{\cX} \exp \left(\frac{f(x)-c(x, y)}{\eps}\right) \frac{d \mu(x)}{dx}dx\right) \\
   &=-\eps \log \left(\int_{\cX} \exp \left(\frac{\tilde{f}(x)-c(x, y)}{\eps}\right)dx\right) + \eps\log(Z)\\
   &= -\eps \log \big(F_{\tilde{f}}(y)\big)+ \eps\log(Z).
\end{align*}

\medskip

We are now ready to prove \eqref{eq:stability_coupling}. For $\rho\in\cP_{\mathsf{ac}}(\cX)$ with Lebesgue density $\frac{d\rho}{dx}$, denote the differential entropy of $\rho$ by $\sH(\rho)\coloneqq -\int_\cX \log\left(\frac{d\rho}{dx}\right)d\rho$. Consider:
\begin{align*}
    &\Gamma\left(\varphi_\star\right)-\Gamma\left(f\right)\\
    &=\int_{\cX \times \cY} c d \pi^\eps_\star-\eps \int_{\cY} \sH\big(\pi^\eps_\star(\cdot| y)\big) d \nu(y)+\eps \sH(\mu)-\int_{\cX} f d \mu-\int_{\cY} f^{c,\eps}d \nu\\
    &=\int_{\cX \times \cY}\big(c(x, y)-\tilde{f}(x)\big) d \pi^\eps_\star(x, y)-\eps \int_{\cY} \sH\big(\pi^\eps_\star(\cdot| y)\big) d \nu(y)+\eps \int_{\cY}\log\big(F_{\tilde f}\big)d\nu\\
     &=-\mspace{-2mu}\eps \mspace{-5mu}\int_{\cX \times \cY} \mspace{-22mu}\frac{\tilde{f}(x)-c(x, y)}{\eps} d \pi^\eps_\star(x, y)\mspace{-2mu}+\mspace{-2mu}\eps \mspace{-5mu}\int_{\cX \times \cY} \mspace{-22mu}\log \big(F_{\tilde{f}}(y)\big) d\pi^\eps_\star(x, y)-\eps \int_{\cY} \sH\big(\pi^\eps_\star(\cdot| y)\big) d\nu(y)\\
     &=-\eps\int_{\cX \times \cY} \log \left(\frac{1}{F_{\tilde{f}}(y)} \exp \left(\frac{\tilde{f}(x)-c(x, y)}{\eps}\right)\right) d\pi^\eps_\star(x, y)-\eps \int_{\cY} \sH\big(\pi^\eps_\star(\cdot| y)\big) d\nu(y)\\
     &=-\eps\mspace{-4mu} \int_{\cX \times \cY} \mspace{-22mu}\log \left(\frac{d\pi^\eps_{f}(x\mid y) }{dx}\right) d \pi^\eps_\star(x, y)\mspace{-2mu}-\mspace{-2mu}\eps \mspace{-5mu}\int_{\cY} \sH\big(\pi^\eps_\star(\cdot| y)\big) d\nu(y)\\
     &=-\eps \int_{\cY}\int_{ \cX} \log \left(\frac{d\pi^\eps_{f}(x| y) }{d x}\right) d \pi^\eps_\star(x| y)d\nu(y)+\eps \int_{\cY}\int_{\cX} \log\left(\frac{d \pi^\eps(x| y)}{d x}\right)d \pi^\eps_\star(x| y) d\nu(y)\\
     &=\eps \int_{\cY}\int_{\cX}\log \left(\frac{d \pi^\eps_\star(x | y)}{d\pi^\eps_{f}(x|y)}\right)d\pi^\eps_\star(x| y)d\nu(y)\\
     &=\eps\int_{\cX \times \cY} \log \left(\frac{d\pi^\eps_\star(x, y)}{d\pi^\eps_{f}(x, y)}\right) d\pi^\eps_\star(x, y)\\
     &=\eps\KL\left(\pi^\eps_\star\middle\|\pi^\eps_{f}\right).
\end{align*}

Recalling that $\hat f_\star$ is a NN that optimizes the NE $\widehat\OT^\eps_{k,a}(X^n,Y^n)$ from \eqref{eq:NE_EOT}, and plugging it into \eqref{eq:stability_coupling} yields $\KL\big(\pi^\eps_\star\big\|\pi^\eps_{\hat{f}_\star}\big)=\eps^{-1}\big[\Gamma_{} \left(\varphi_\star\right)-\Gamma(\hat{f}_\star)\big]$. Thus, to prove the KL divergence bound from Theorem \ref{thm:Bound_NeuralCoupling}, it suffices to control the gap between the $\Gamma$ functionals on the RHS above.

\medskip
Write $f_{\star}$ for a NN that maximizes the population-level neural EOT cost $\OT_{k,a}^\eps(\mu,\nu)$ (see \eqref{eq:neural_pop_EOT}). Define $\widehat{\Gamma}(f): = \frac 1n\sum_{i=1}^n\big[f(X_i)+f^{c, \eps}(Y_i)\big]$ for the optimization objective in the problem $\widehat\OT_{k,a}^\eps(X^n,Y^n)$ (see \eqref{eq:NE_EOT}), and note that $\hat{f}_{\star}$ is a maximizer of $\widehat{\Gamma}$. We now have
\begin{equation}
\begin{aligned}
&\eps\KL\left(\pi_\star^\eps\middle\|\pi^\eps_{\hat{f}_{\star}}\right)\\
& = \Gamma\left(\varphi_{\star}\right)-\Gamma \left(f_{\star}\right)+\Gamma\left(f_\star\right)-\widehat{\Gamma} \left(\hat{f}_\star\right)+\widehat{\Gamma} \left(\hat{f}_\star\right)-\Gamma\left(\hat{f}_\star\right)\\
&=\underbrace{\OT^{\eps}(\mu, \nu)\mspace{-2mu}-\mspace{-2mu}\OT^{\eps}_{k,a}(\mu, \nu)}_{(\mathrm{I})} + \underbrace{\OT^{\eps}_{ k,a}(\mu, \nu)\mspace{-2mu}-\mspace{-2mu}\widehat\OT^{\eps}_{k,a}\left(X^n, Y^n\right)}_{(\mathrm{II})}+\underbrace{\widehat{\Gamma} \left(\hat{f}_\star\right)\mspace{-2mu}-\mspace{-2mu}\Gamma \left(\hat{f}_\star\right)}_{(\mathrm{III})}.
\end{aligned}
\end{equation}

Setting $a =\bar{c}_{b,d}$ as in the proof of Theorem \ref{thm:Bound_EGW1} (see \eqref{eq:constant_b} and \eqref{eq:constant_c_bar}) and taking an expectation (over the data) on both sides, Terms $(\mathrm{I})$ and $(\mathrm{II})$ are controlled, respectively, by the approximation error and empirical estimation error from \eqref{eq:EGW_approx_error} and \eqref{eq:EGW_est_error}. For Term $(\mathrm{III})$, consider
\begin{align*}
    \mathbb{E}\left[\widehat{\Gamma} \left(\hat{f}_\star\right)-\Gamma \left(\hat{f}_\star\right)\right]
    &\leq \mathbb{E}\left[\sup_{f \in \cF_{k, d}(\bar{c}_{b,d})} \left|\Gamma(f)-\widehat{\Gamma}(f)\right|\right]\\
    &\leq n^{-\frac{1}{2}}\mathbb{E}\left[\sup _{f \in \cF_{k, d}(\bar{c}_{b,d})}\mspace{-10mu} n^{-\frac{1}{2}}\left|\sum_{i=1}^n\big(f\left(X_i\right)-\mathbb{E}_\mu[f]\big)\right|\right]\mspace{-2mu}\mspace{-2mu}\\
    &\quad\quad\quad+\mspace{-2mu}n^{-\frac{1}{2}}\mathbb{E}\left[\sup _{f \in \mathcal{F}^{c, \eps}} \mspace{-3mu}n^{-\frac{1}{2}}\left|\sum_{j=1}^n\big(f\left(Y_j\right)-\mathbb{E}_\nu[f]\big)\right|\right].
\end{align*}
where the last step follows similarly to \eqref{eq:emp_error_decomposition}. Notably, the RHS above is also bounded  by the estimation error bound from \eqref{eq:EGW_est_error}. Combining the above completes the proof.
\qed

\section{Proofs of Technical Lemmas}\label{appen:proofs_lemma}
\subsection{Proof of Lemma \ref{lemma:regularity_EOT}}
\label{appen:proof:lemma:regularity_EOT}
    The existence of optimal potentials follows by standard EOT arguments \cite[Lemma 1]{goldfeld2022limit}. Recall that EOT potentials are unique up to additive constants. Thus, let $\left(\varphi_0, \psi_0\right) \in L^1(\mu) \times L^1(\nu)$ be optimal EOT potentials for the cost $c$, solving dual formulation \eqref{eq:EOT_dual}, and we can assume without loss of generality that $\int \varphi_0 d\mu=\int \psi_0 d\nu=\frac{1}{2} \OT^{\mspace{1mu}\eps}(\mu, \nu)$.

    Recall that the optimal potentials satisfies the Schr\"odinger system from \eqref{eq:Schrodinger system}. Define new functions $\varphi$ and $\psi$ as
$$
\begin{array}{ll}
\varphi(x)\coloneqq-\eps\log \int_{\cY} \exp\left(\frac{\psi_0(y)-c(x, y)}{\eps}\right) d\nu(y), & x \in \cX \\
\psi(y)\coloneqq \varphi^{c ,\eps}(y), & y \in \cY.
\end{array}
$$
These integrals are clearly well-defined as the integrands are everywhere positive on $\cX$ and $\cY$, and $\varphi_0, \psi_0$ are defined on the supports of $\mu, \nu$ respectively. Now We show that $\varphi,\psi$ are pointwise finite. For the upper bound, by Jensen's inequality, we have
$$
\varphi(x) \leq \int_{\cY} \frac{1}{2}\|x-y\|^2-\psi_0(y) d \nu(y) \leq 2d,
$$
the second inequality follows from $\int \psi_0 d\nu=\frac{1}{2} \OT^{\mspace{1mu}\eps}(\mu, \nu) \geq 0$. The upper bound holds similarly for $\psi$ on $\cY$ and $\psi_0$ on the support of $\nu$. For lower bound, we have
\begin{align*}
    -\varphi(x) &\leq \eps\log\int_{\cY} \exp\left(\frac{2 d}{\eps}\right) d\nu(y)=2d.
\end{align*}
Note that $\varphi$ is defined on $\cX$, with pointwise bounds proven above. By Jensen's inquality,
$$
\begin{aligned}
\int_{\cX}\left(\varphi_0-\varphi\right) d \mu+\int_{\cY}\left(\psi_0-\psi\right) d\nu
& \leq \eps\log \int_{\cX} \exp\left(\frac{\varphi_0-\varphi}{\eps}\right) d\mu+\log \int_{\cY} \exp\left(\frac{\psi_0-\psi}{\eps}\right) d \nu \\
& =\eps\log \int_{\cX \times \cY} \exp\left(\frac{\varphi_0(x)+\psi_0(y)-c(x, y)}{\eps}\right) d \mu \otimes \nu \\
& \quad\quad+\eps\log \int_{\cX \times \cY} \exp\left(\frac{\varphi(x)+\psi_0(y)-c(x, y)}{\eps}\right) d \mu \otimes \nu\\
& =0 .
\end{aligned}
$$
Since $\left(\varphi_0, \psi_0\right)$ maximizes \eqref{eq:EOT_dual}, so does $(\varphi, \psi)$ and thus they are also optimal potentials. Therefore, $\varphi$ solves semi-dual formulation \eqref{eq:EOT_semidual}. By the strict concavity of the logarithm function we further conclude that $\varphi=\varphi_0$ $ \mu$-a.s and $\psi=\psi_0$ $\nu$-a.s.

The differentiability of $(\varphi, \psi)$ is clear from their definition. For any multi-index $\alpha$, the multivariate Faa di Bruno formula (see \cite[Corollary 2.10]{constantine1996multivariate}) implies
\begin{equation}
\label{eq:constant_C_s_1}
\begin{aligned}
       -D^\alpha \varphi(x)=\eps\mspace{-3mu}\sum_{r=1}^{|\alpha|}\mspace{-2mu}\sum_{p(\alpha, r)}\mspace{-5mu}\frac{\alpha !(r-1) !(-1)^{r-1}}{\prod_{j=1}^{|\alpha|}\left(k_{j} !\right)\left(\beta_{j} !\right)^{k_j}}\mspace{-2mu}\prod_{j=1}^{|\alpha|}\mspace{-2mu}\left(\mspace{-2mu}\frac{D^{\beta_j}\mspace{-7mu} \int \exp(\frac{\psi_0(y)-c(x, y)}{\eps}) d \nu(y)}{\int \exp(\frac{\psi_0(y)-c(x, y)}{\eps}) d \nu(y)}\mspace{-2mu}\right)^{k_j}, 
\end{aligned}
\end{equation}
where $p(\alpha, r)$ is the collection of all tuples $\left(k_1, \cdots, k_{|\alpha|} ; \beta_1, \cdots, \beta_{|\alpha|}\right) \in \mathbb{N}^{|\alpha|} \times \mathbb{N}^{d \times|\alpha|}$ satisfying $\sum_{i=1}^{|\alpha|} k_i=r, \sum_{i=1}^{|\alpha|} k_i \beta_i=\alpha$, and for which there exists $s \in\{1, \ldots,|\alpha|\}$ such that $k_i=0$ and $\beta_i=0$ for all $i=1, \ldots|\alpha|-s, k_i>0$ for all $i=|\alpha|-s+1, \ldots,|\alpha|$, and $0 \prec \beta_{|\alpha|-s+1} \prec \cdots \prec \beta_{|\alpha|}$. For a detailed discussion of this set including the linear order $\prec$, please refer to \cite{constantine1996multivariate}. For the current proof we only use the fact that the number of elements in this set solely depends on $|\alpha|$ and $r$. Given the above, it clearly suffices to bound $\left|D^{\beta_j} \int \exp\left(\psi_0(y)-c(x, y)\right) d\nu(y)\right|$. First, we apply the same formula to $D^{\beta_j} e^{-c(x, y)/\eps}$ and obtain
\begin{equation}
\label{eq:constant_C_s_2}
\begin{aligned}
      D^{\beta_j} e^{\frac{-c(x, y)}{\eps}}=\sum_{r^{\prime}=1}^{\left|\beta_j\right|}(\frac{-1}{2\eps})^{r^\prime}\mspace{-10mu}\sum_{p\left(\beta_j, r^{\prime}\right)} \mspace{-8mu}l(\beta_j,r^\prime,\boldsymbol{k}^\prime,\boldsymbol{\eta}) e^{\frac{-c(x, y)}{\eps}} \prod_{i=1}^{\left|\beta_j\right|}\mspace{-2mu}\left(D^{\eta_i}\|x\mspace{-1mu}-\mspace{-1mu}y\|^2\right)^{k_i^{\prime}} ,
\end{aligned}
\end{equation}
where $p\left(\beta_j, r^{\prime}\right)$ defined similarly to the above. Observe that
\begin{align*}
   \left|D^{\eta_i}\left(\|x-y\|^2\right)\right| &\leq 4\left(1+\|x-y\|\right)\leq 4(1+2\sqrt{d}),
\end{align*}
where the first inequality follows from proof of Lemma 3 in \cite{zhang2022gromov}. Consequently, for $0<\eps< 1$, we have
$$
    \left|\frac{D^{\beta_j} \int \exp(\frac{\psi_0(y)-c(x, y)}{\eps}) d \nu(y)}{\int \exp(\frac{\psi_0(y)-c(x, y)}{\eps}) d \nu(y)}\right| \leq C_{\beta_j}\eps^{-|\beta_j|}(1+2\sqrt{d})^{|\beta_j|},
$$
and for $\eps \geq1$,
$$
    \left|\frac{D^{\beta_j} \int\exp(\frac{\psi_0(y)-c(x, y)}{\eps}) d \nu(y)}{\int \exp(\frac{\psi_0(y)-c(x, y)}{\eps}) d \nu(y)}\right| \leq C_{\beta_j}\eps^{-1}(1+2\sqrt{d})^{|\beta_j|},
$$
Plugging back, we obtain 
$$
\left|D^\alpha \varphi(x)\right| \leq C_{|\alpha|}(1+\eps^{1-|\alpha|})(1+2\sqrt{d})^{|\alpha|}.
$$
Analogous bound holds for $\psi$.
\qed

\subsection{Proof of Lemma \ref{lemma:neural_approx_reduction}}
\label{appen:proof_neural_approx_reduction}

For any $f\in \cF_{k,a}(a)$, we know that $\|f\|_{\infty,\cX} \leq 3 a(\|\cX\|+1)=6a$, so NNs are uniformly bounded. This implies that $\OT^{\eps}(\mu,\nu) \geq \OT^{\eps}_{k,a}(\mu,\nu)$. Since $\varphi$ satisfies \eqref{eq:regularity_EOT}, it's uniformly bounded on $\cX$. Then, the following holds:
\begin{align*}
\left|\OT^{\eps}(\mu,\nu)-\OT^{\eps}_{k,a}(\mu,\nu)\right|  =&\OT^{\eps}(\mu,\nu)-\OT^{\eps}_{k,a}(\mu,\nu) \\
\leq &\mathbb{E}_{\mu} |\varphi-f|+\mathbb{E}_{\nu} |\varphi^{c,\eps}-f^{c,\eps}|\\
\leq & 2\|\varphi-f\|_{\infty,\cX}.
\end{align*}  
The last inequality holds by an observation that, 
\begin{equation}
\label{eq:lipshcitz_c_tranform}
    |\varphi^{c,\eps}(y)-f^{c,\eps}(y)| \leq \left\|\varphi-f\right\|_{\infty, \cX},\quad \forall y\in \cY.
\end{equation}
Indeed, note that for any $y \in \cY$,
\begin{align*}
    \varphi^{c,\eps}(y)-f^{c,\eps}(y)
    &= -\eps \log \left(\frac{\int \exp \left(\frac{\varphi(x)-c(x, y)}{\eps}\right) d \mu(x)}{\int \exp \left(\frac{f(x)-c(x, y)}{\eps}\right) d \mu(x)}\right)\\
    & = -\eps \log \left(\frac{\int \exp\left(\frac{\varphi(x)-f(x)}{\eps}\right)\exp \left(\frac{f(x)-c(x, y)}{\eps}\right) d \mu(x)}{\int\exp \left(\frac{f(x)-c(x, y)}{\eps}\right) d \mu(x)}\right)\\
    & \geq -\eps \log \left(\frac{\int \exp\left(\frac{\|\varphi-f\|_{\infty,\cX}}{\eps}\right)\exp \left(\frac{f(x)-c(x, y)}{\eps}\right) d \mu(x)}{\int\exp \left(\frac{f(x)-c(x, y)}{\eps}\right) d \mu(x)}\right)\\
    & = -\left\|\varphi-f\right\|_{\infty, \cX}.
\end{align*}
Similarly, we can have that
\begin{align*}
    \varphi^{c,\eps}(y)-f^{c,\eps}(y) 
    & \leq -\eps \log \left(\frac{\int\exp\left(\frac{-\|\varphi-f\|_{\infty,\cX}}{\eps}\right)\exp \left(\frac{f(x)-c(x, y)}{\eps}\right) d \mu(x)}{\int \exp \left(\frac{f(x)-c(x, y)}{\eps}\right) d \mu(x)}\right)\\
    & = \left\|\varphi-f\right\|_{\infty, \cX}.
\end{align*}
\qed

\end{document}